\documentclass[12pt,letter]{article} 
\usepackage[T2A]{fontenc}
\usepackage[cp866]{inputenc}
\usepackage[russian,english]{babel}
\usepackage[tbtags]{amsmath}
\usepackage{amsfonts,amssymb,mathrsfs,amscd,comment}
\usepackage{amsthm,graphicx} 
\usepackage{hyperref}
\def\NoBlackBoxes{\overfullrule0pt}

\NoBlackBoxes
\theoremstyle{plain}
\newtheorem{theorem}{Theorem}

\theoremstyle{definition}

\newtheorem{remark}{Remark}

\theoremstyle{main}
\def\bad{\spaceskip=0.33emplus0.6emminus0.15em\immediate\write5{\string\bad}}

\theoremstyle{plain}
\newtheorem{proposition}{Proposition}
\let\myh\widehat
\let\myo\overline

\let\pfi\varphi

\def\RR{\mathbb R}
\def\CC{\mathbb C}

\def\NN{\mathbb N}

\def\const{\operatorname{const}}
\def\<{\left\langle}
\def\>{\right\rangle}
\def\({\left(}
\def\){\right)}
\def\[{\left[}
\def\]{\right]}
\NoBlackBoxes

\def\bad{\spaceskip=0.33emplus0.6emminus0.15em\immediate\write5{\string\bad}}

\def\NN{\mathbb N}
\def\RR{\mathbb R}
\def\CC{\mathbb C}
\def\PP{\mathbb P}
\def\FF{\mathbf F}
\def\zz{\mathbf z}
\def\RS{\mathfrak R}
\def\vv{V}
\def\bmu{{\boldsymbol\mu}}
\def\maa{\mathbf a}
\def\mytt{\mathbf t}
\def\blambda{{\boldsymbol\lambda}}
\def\bGamma{\boldsymbol\Gamma}
\def\hp{\vphantom{Hp}}
\def\dd{\mathrm{dd}}

\let\myh\widehat

\let\myo\overline
\def\const{\operatorname{const}}
\def\Nik{\operatorname{Nik}}
\let\pfi\varphi
\def\({\left(}
\def\){\right)}
\def\[{\left[}
\def\]{\right]}
\def\<{\left\langle}
\def\>{\right\rangle}

\begin{document}

\selectlanguage{english}


\author{Sergey~P.~Suetin}

\title{On the equivalence of the scalar and vector equilibrium problems for a~pair of
functions forming a~Nikishin system}




\maketitle

\markright{}

\begin{abstract}
We prove the equivalence of the vector and scalar equilibrium problems which
arise naturally in the study of the limit zeros distribution of type~I
Hermite--Pad\'e polynomials for a~pair of functions forming a~Nikishin system.

Bibliography: 22 titles.
\end{abstract}



\footnotetext[1]{This work is supported by the Russian Science Foundation
(grant no.~19-11-00316).}

\section{Introduction and statement of the problem}\label{s1}
\subsection{}\label{s1s1}
The purpose of the present paper is to further develop a~new approach
to the study of extremal and equilibrium problems that appear naturally
when examining the limit zeros distribution of Hermite--Pad\'e
polynomials. The crux of this approach, which was proposed by the
author of the present paper in~\cite{Sue18} (see also~\cite{Sue18b}),
is to consider, instead of the traditional vector equilibrium problem
on the Riemann sphere~$\myh{\CC}$, the scalar problem (but already on
a~Riemann surface). Here we give some arguments supporting the
naturality and expedience of this alternative scalar approach. Namely,
without having recourse to the problem on the limit zeros distribution
of Hermite--Pad\'e polynomials
 and based only on the potential theory on a~compact Riemann surface,
 we prove that in~\eqref{19}, as well as in the vector equilibrium problem~\eqref{5},
 the equality takes place on the whole of the compact set~$\FF$ (which implies that
$S(\lambda_{\FF})=\FF$). Moreover, we show here that the vector and the scalar problems
are in a~sense equivalent (see Theorem~\ref{th1} below).
In subsequent studies, we are also planning to prove the existence of the limit
zeros distribution of type~II Hermite--Pad\'e polynomials based only on the appropriate scalar
equilibrium problem in potential theory posed on a~Riemann surface.

Note that Stahl \cite{Sta87} and~\cite{Sta88} proposed a~certain approach to the above problems
by employing the machinery of the potential theory on a~compact Riemann surface.
However, his approach has not been worked out.
The approach put forward by the author of the present paper in~\cite{Sue18} is different from Stahl's approach. In particular,
as distinct from the author's papers \cite{Sue18} and~\cite{Sue18b}, Stahl \cite{Sta87}
and~\cite{Sta88} has never considered
extremal problems of the potential theory or equilibrium problems, even though he worked with
potentials on a~compact Riemann surface.

Following ~\cite{Sue18} (see also~\cite{Sue18b}), consider
\begin{equation}
f_1(z):=\frac1{(z^2-1)^{1/2}},\quad
f_2(z):=\frac1\pi\int_{-1}^1\frac{h(x)}{z-x}\frac{dx}{\sqrt{1-x^2}},
\quad z\in D:=\myh{\CC}\setminus E,
\label{1}
\end{equation}
where $E:=[-1,1]$ and the branch of the function $(\,{\cdot}\,)^{1/2}$ is chosen such that $(z^2-1)^{1/2}/z\to1$ as
$z\to\infty$; by $\sqrt{1-x^2}$, $x\in(-1,1)$, we mean the positive square root: $\sqrt{a^2}=a$ for $a\geq0$.
Here and in what follows we assume that $h=\myh{\sigma}$ in~\eqref{1} is
a~Markov function
supported on a~compact set $F\subset\RR\setminus E$; i.e.,
\begin{equation}
h(z)=\myh{\sigma}(z):=\int_{F}\frac{d\sigma(t)}{z-t},\quad z\in\myh\CC\setminus
F,
\label{2}
\end{equation}
where
$$
F:=\bigsqcup_{j=1}^m[c_j,d_j]\subset\RR\setminus E,
$$
$c_1<d_1<\dots< c_m<d_m$, $\sigma$ is a~positive Borel measure with support $S(\sigma)$ on~$F$ and
such that $S(\sigma)=F$ and
$\sigma':=d\sigma/dx>0$ almost everywhere on~$F$ (see~\cite{Sue18}).
Throughout we will use these notation and conventions.

For $f_1$ we have the representation
\begin{equation}
f_1(z)=\frac1\pi\int_{-1}^1\frac1{z-x}\frac{dx}{\sqrt{1-x^2}},\quad z\in
D,
\label{3}
\end{equation}
and hence using \eqref{1},~\eqref{2} and~\eqref{3} we see that
$\Delta f_2(x)/\Delta f_1(x)=\myh\sigma(x)$ for $x\in(-1,1)$,
where $\Delta f_j(x)$ denotes the difference of the limit values of the function
$f_j$, $j=1,2$, taken from the upper and lower half-planes, respectively.
It follows that the pair of functions $(f_1,f_2)$ forms a~Nikishin system
(for more details on this concept, see \cite{Nik86},~\cite{NiSo88},
and also~\cite{ApBoYa17}, \cite{BaGeLo18},~\cite{LoVa18}, as well as
the references given therein).
In~\cite{Sue18c} an example of a~multivalued analytic function~$f$ is given
such that the pair $f,f^2$ forms a~Nikishin system (note that in~\cite{Sue18c} the concept of a~Nikishin
system is a~little bit more general than that given by E.~M.~Nikishin himself).
There exist classes of multivalued analytic functions~$f$ such that
the pair of functions $f$,~$f^2$ can be naturally looked upon as a~{\it complex} Nikishin system
(see \cite{RaSu13}, \cite{MaRaSu16},~\cite{Sue18c}). In connection with the new approach of~\cite{Sue18d}
to the problem of efficient continuation of a~given germ of a~multivalued analytic function,
this fact seems to be one of the main impetus for the study of equilibrium problems pertaining to complex Nikishin systems.

Given an arbitrary $n\in\NN$, we denote by $\PP_n$ the set of all polynomials of degree
$\leq{n}$ with complex coefficients.
For an arbitrary polynomial $Q\in\PP_n^{*}:=\PP_n\setminus0$, we let
$\chi(Q)$ denote
 the measure counting the zeros (with multiplicities) of the polynomial~$Q$,
$$
\chi(Q):=\sum_{\zeta:Q(\zeta)=0}\delta_\zeta;
$$
$\delta_\zeta$ is the unit measure concentrated at the point $\zeta$ (the Dirac delta-measure).

For a tuple of three functions $[1,f_1,f_2]$, where $f_1$ and $f_2$ are given by
representations \eqref{1}, and for an arbitrary $n\in\NN$, the type~I Hermite--Pad\'e polynomials
$Q_{n,0},Q_{n,1},Q_{n,2}\in\PP^{*}_n$ are defined from the relation
\begin{equation}
R_n(z):=(Q_{n,0}+Q_{n,1}f_1+Q_{n,2}f_2)(z)=O\(\frac1{z^{2n+2}}\),
\quad z\to\infty,
\label{4}
\end{equation}
in the standard way.
It is well known that such polynomials always exist, but they are not uniquely
specified by~\eqref{4}; for more details, see \cite{Nut84}, \cite[Ch.~4, \S~1]{NiSo88},
\cite{ApBoYa17}, and~\cite{BaGeLo18}.
The solution of the problem on the limit (as $n\to\infty$) zeros distribution
of the Hermite--Pad\'e polynomials $Q_{n,j}$ for a~pair of functions $f_1,f_2$
defined by~\eqref{1} can be obtained from the results of E.~M.~Nikishin~\cite{Nik86}
(see also~\cite{NiSo88}). In~\cite{Nik86} this problem was solved by E.~M.~Nikishin
(in a~much more general setting than the one considered in the present paper)
in terms of the vector equilibrium problem with a~$2\times2$-interaction matrix (which is now called
a~{\it Nikishin matrix}) on the basis of the general {\it vector} approach,
which was first proposed by A.~A.~Gonchar and
E.~E.~Rakhmanov~\cite{GoRa81}. For further advances in this vector approach, see \cite{Apt08}, \cite{Rak12}, \cite{RaSu13},
and~\cite{ApBoYa17}.
At the same time, the vector approach faces certain difficulties in the solution of problems
involving complex Nikishin systems of the form $f,f^2$, where, for example, $f$~is a~multivalued function
of Laguerre class (see~\cite{MaRaSu16}). In~\cite{Sue18} (see also~\cite{Sue18b})
a~new approach to the solution of the problem on the limit zeros distribution
of the Hermite--Pad\'e polynomials
of type~I was proposed. This approach is based on the solution of the scalar equilibrium problem,
but this problem is posed not on the Riemann sphere~$\myh{\CC}$, but rather
on a~two-sheeted Riemann surface (briefly, RS) of the function $w^2=z^2-1$.
This approach proved instrumental in delivering, by a~different method, the results
established earlier in~\cite{Nik86} and~\cite{RaSu13};
it has become possible to derive some new results not amenable to the vector approach machinery
(see~\cite{Sue18b}). The purpose of the present paper is to prove, for a~given pair
of compact sets $E=[-1,1]$ and $F=\bigsqcup_{j=1}^m[c_j,d_j]$, the equivalence of the traditional
vector equilibrium problem related to a~Nikishin system (see \S\,\ref{s1s2}
below) and the scalar problem \cite[(1.17)]{Sue18}, which was posed in~\cite{Sue18}.
Furthermore, it will be shown that $S(\lambda_\FF)=\FF$ in \cite[(1.17)]{Sue18}; thereby
we prove that in the equilibrium relations (see~\eqref{19}) there is an equality sign
on the whole of the compact set $\FF$. Of course, the equivalence of the
equilibrium problems also follows from the fact that
both the solution of one problem and the solution of the other problem
characterize the limit zeros distribution of the same
Hermite--Pad\'e polynomials of type~I. Here we prove the equivalence directly in the terms
related to these equilibrium problems and without recourse to Hermite--Pad\'e polynomials.

\subsection{}\label{s1s2}
We recall some well-known facts pertaining to the vector equilibrium problem appearing in the solution
of the problem on the limit zeros distribution of Hermite--Pad\'e polynomials for a~pair
of functions forming a~Nikishin system (for more details, see \cite{Nik86}, \cite{NiSo88},~\cite{BaGeLo18}).
In accordance to the vector approach, which dates back to A.~A.~Gonchar and E.~A.~Rakhmanov \cite{GoRa81},
the answer to the problem on the limit zeros distribution of Hermite--Pad\'e polynomials
is given precisely in terms related to the unique vector measure in which this equilibrium problem is solved.

For an arbitrary (positive Borel) measure $\mu$, $S(\mu)\subset\CC$, by
$$
U^\mu(z):=\int\log\frac1{|z-t|}\,d\mu(t),\quad
z\in\CC\setminus S(\mu),
$$
we denote the logarithmic potential of the measure~$\mu$.

Let, as before, $E=[-1,1]$ and $F$~be a~compact set consisting of a~finite number of
closed intervals lying on the real line, $E\cap F=\varnothing$.
We let $M_1(E)$ and $M_1(F)$ denote, respectively, the space of all unit
measures supported on $E$ and~$F$. We also denote by $M^\circ_1(E)$ and $ M^\circ_1(F)$
the subspaces of the spaces $M_1(E)$ and $M_1(F)$, respectively, with finite energy (with respect to the logarithmic kernel). Let
$$
M_{\Nik}=
\begin{pmatrix} 4&-1\\-1&1\end{pmatrix}
$$
be the Nikishin $2\times2$-matrix and let
\begin{equation}
\left\{
\begin{aligned}
4U^{\lambda_1}(x)-U^{\lambda_2}(x)\equiv w_1 &=\const,\quad x\in E,\\
-U^{\lambda_1}(t)+U^{\lambda_2}(t)\equiv w_2 &=\const,\quad t\in F,
\end{aligned}\right.
\label{5}
\end{equation}
be the corresponding vector $2\times2$-equilibrium problem with respect to the
measures $\lambda_1\in M_1^\circ(E)$ and $\lambda_2\in M_1^\circ(F)$.
It is well known that the vector measure $\vec{\lambda}:=(\lambda_1,\lambda_2)$, which is
the solution of the equilibrium problem~\eqref{5}, exists and is unique.
Moreover, $S(\lambda_1)=E$ and $S(\lambda_2)=F$.
The second of~\eqref{5} implies that
$\lambda_2=\beta_F(\lambda_1)$ is the balayage of the measure $\lambda_1\in M_1(E)$ from the domain
$\Omega:=\myh{\CC}\setminus F$ onto its boundary $\partial\Omega=F$, and hence
problem~\eqref{5} is equivalent to the following equilibrium problem (see \cite{RaSu13}, \cite{BuSu15},~\cite{MaRaSu16}):
\begin{equation}
3U^{\lambda_1}(x)+G^{\lambda_1}_F(x)\equiv\const,\quad x\in E;
\label{6}
\end{equation}
here
\begin{equation}
G^{\mu}_F(z):=\int g_F(z,\zeta)\,d\mu(\zeta)
\label{7}
\end{equation}
is the Green potential of the measure~$\mu$, $S(\mu)\subset\CC$, and $g_F(z,\zeta)$ is the Green function
for the domain~$\Omega$ with singularity at the point $z=\zeta$.
From a~unique measure~$\lambda_1$ satisfying \eqref{6} the
measure $\lambda_2$ is recovered in a~unique way. Namely,
$\lambda_2=\beta_F(\lambda_1)$.
Moreover (see \cite[Ch.~5, \S~7, (7.13)]{NiSo88} and also~\cite{ApBoYa17}),
$\deg{Q_{n,2}}=n$, all zeros of the polynomial $Q_{n,2}$ lie in the convex hull~$\myh{F}$ of the compact set~$F$,
and
\begin{equation}
\frac1n\chi(Q_{n,2})\overset{*}\longrightarrow\lambda_2,\quad n\to\infty
\label{10}
\end{equation}
(see \cite[Ch.~5, \S~7, Theorems 7.1 and 7.4]{NiSo88});
here and in what follows ``$\overset{*}\longrightarrow$'' denotes the weak-$*$
convergence in the space of measures.

So, problem \eqref{6} is an equilibrium problem for one measure
$\lambda_1$, rather than for two measures $\lambda_1$ and $\lambda_2$ (as problem~\eqref{5});
recall that the measure $\lambda_2$ is now uniquely defined from the measure~$\lambda_1$.
Nevertheless, the equilibrium problem~\eqref{6} still should be looked upon as a~vector
problem, because its statement depends on both compact sets:
the compact set~$E$ and the compact set~$F$.
This is the first reason why the attempts to extend this problem
to the complex\footnote{An equilibrium problem will be called complex if
not all branch points of the function~$h$ lie on the real line; see \cite{RaSu13},~\cite{Sue18b}.}
setting involve considerable difficulties -- under this approach
one has in fact at first to find two $S$-curves: the one replacing the closed interval~$E$ and
the other one replacing the compact set~$F$
(see \cite{Nut84},~\cite{Rak12},~\cite{Apt08},~\cite{RaSu13}).
In contrast, the advantage of the scalar approach, which was introduced in~\cite{Sue18}, is that
in the complex setting it leads to the problem of finding a~single
$S$-curve; but this curve should lie on the RS
(see \cite{Sue18}, \cite{Sue18b},~\cite{Sue18c}).
The corresponding equilibrium problem is now phrased in terms of some potential
on a~compact RS,
and in general, in the presence of a~harmonic external field (see \cite{RaSu13},~\cite{Chi18},~\cite{Chi19}).

The scalar approach discussed here was first proposed by the author~\cite{Sue18}
and applied to the solution of the
problem on the limit zeros distribution of the Hermite--Pad\'e polynomials
$Q_{n,2}$ which are defined from~\eqref{4} and constructed for the pair
of functions $f_1,f_2$ defined by~\eqref{1}.
We note once more that in this specific setting the crux of the scalar approach
is that the equilibrium problem is posed on the two-sheeted RS $\RS_2=\RS_2(w)$
of the function $w^2=z^2-1$ and in potential-specific terms on a~compact RS
(cf.~\cite{ApBoYa17},~\cite{LoVa18}).

\begin{remark}\label{rem1}
As was pointed out in \cite{Sue18}, in the case considered here the solution to the problem
on the limit zeros distribution of Hermite--Pad\'e polynomials was obtained already
by E.~M.~Nikishin~\cite{Nik86} in the framework of the traditional approach. So,
both in~\cite{Sue18} and in the present paper we speak about the proof of the equivalence of the new
scalar approach and the traditional vector approach on an example of the previously solved problems.
The advantages of the scalar approach are manifested in the solution of complex problems
not amenable to the machinery of the conventional vector approach; see~\cite{Sue18b}.
\end{remark}

\subsection{}\label{s1s3}
We require the following notation and definitions from~\cite{Sue18}.

Given $z\in D=\myh{\CC}\setminus E$, we denote by
\begin{equation}
\pfi(z):=z+(z^2-1)^{1/2}
\label{11}
\end{equation}
the inverse Zhukovskii function (recall that everywhere in the present
paper we choose a~branch of the function $(\,{\cdot}\,)^{1/2}$ such
that $(z^2-1)^{1/2}/z\to1$ as $z\to\infty$). The function~$\pfi$ is
a~single-valued meromorphic function in the domain~$D$. We assume that
a~point $\zz\in\RS_2$ has the form $\zz=(z,w)$. Let
$\pi_2\colon\RS_2\to\myh\CC$ be the two-sheeted covering ($\pi_2$~is
the canonical projection), $\pi_2(\zz)=z$. The function $\pfi$ is
defined on the RS~$\RS_2$ by the equality $\pfi(\zz):=z+w$.

We decompose the RS $\RS_2$ into two open sheets $\RS_2^{(0)}$ (the
zero\footnote{As usual,
we identify the sheet $\RS_2^{(0)}$ of the RS $\RS_2$ with the ``physical''
domain $D=\myh{\CC}\setminus E$ of the Riemann sphere.}
sheet) and $\RS_2^{(1)}$ (the first sheet) as follows:
$z^{(0)}:=(z,(z^2-1)^{1/2})\in\RS_2^{(0)}$,
$z^{(1)}:=(z,-(z^2-1)^{1/2})\in\RS_2^{(1)}$.
Setting
$g_2(\zz):=-\log|\pfi(\zz)|=\log|z-w|$, we have $g_2(z^{(0)})=-\log|z|+O(1)$,
$g_2(z^{(1)})=\log|z|+O(1)$, $z\to\infty$, and $g_2(z^{(0)})<g_2(z^{(1)})$.
So, the above partition of the RS $\RS_2$ into sheets is a~Nuttall partition
(see \cite[\S\,3]{Nut84}).
Moreover, $\pi_2(\RS_2^{(0)})=\pi_2(\RS_2^{(1)})=D$.

We set $\vv(\zz):=-\log|\pfi(\zz)|=\log|z-w|$ for $\zz\in\RS_2$; in what follows,
the function $\vv(\zz)$ will play the role of an external field\footnote{More precisely,
$\vv(\zz)$ will play the role of a~potential of the external field.} in the
equilibrium problem considered here.
Let $\FF=F^{(1)}\subset\RS_2$ be some compact set from the first sheet
$\RS^{(1)}_2$ of the RS $\RS_2$
and such that $\pi_2(\FF)=F$. By $M_1(\FF)$ we denote the space of all
unit (positive Borel) measures supported on~$\FF$.
Following~\cite{Sue18}, given an arbitrary measure $\bmu\in M_1(\FF)$,
we introduce the function $P^{\bmu}(\zz)$ (the ``potential'' of the measure~$\bmu$;
see Remark~\ref{rem3}
below) of a~point $\zz\in\RS_2$ by
\begin{equation}
P^{\bmu}(\zz):=\int_{\FF}\log\frac{\left|1-1/\bigl(\pfi(\zz)\pfi(t^{(1)})\bigr)\right|}
{|z-t|^2}
\,d\bmu(t^{(1)}),\quad \zz\in\RS_2\setminus(F^{(0)}\cup F^{(1)}),
\label{12}
\end{equation}
and define the corresponding energy of the measure $\bmu$ (cf.~\cite{Chi18} and~\cite{Chi19})
\begin{equation}
J(\bmu):=
\iint_{\FF\times\FF}\log\frac{\left|1-1/\bigl(\pfi(\zz)\pfi(\mytt)\bigr)\right|}
{|z-t|^2}\,d\bmu(z^{(1)})\,d\bmu(t^{(1)})
=\int_{\FF} P^\bmu(\zz)\,d\bmu(\zz)
\label{14}
\end{equation}
with respect to the kernel
$$
\log\frac{\left|1-1/\bigl(\pfi(\zz)\pfi(\mytt)\bigr)\right|}{|z-t|^2}.
$$
We also consider the energy of the measure~$\bmu$ in the external field $\vv$:
\begin{align}
J_\vv(\bmu):&=
\iint_{\FF\times\FF}\biggl\{
\log\frac{\left|1-1/\bigl(\pfi(\zz)\pfi(\mytt)\bigr)\right|}
{|z-t|^2}
+\vv(\zz)+\vv(\mytt)\biggr\}\,d\bmu(z^{(1)})\,d\bmu(t^{(1)})\notag\\
&=\int_{\FF} P^\bmu(\zz)\,d\bmu(\zz)+2\int_{\FF}\vv(\zz)\,d\bmu(\zz).
\label{15}
\end{align}

By $M_1^\circ(\FF)$ we denote the set of all measures $\bmu\in M_1(\FF)$
with finite\footnote{By~\eqref{14}, this set corresponds with the set of all probability measures
with support on~$F$ and having finite energy with respect to the logarithmic kernel $-\log|z-t|$.}
energy $J(\bmu)$. In what follows, we identity the measure $\bmu\in M_1(\FF)$ and
the measure $\mu=\pi_2(\bmu)\in M_1(F)$,
where $\pi_2(\bmu)(e)=\bmu(e^{(1)})$ for any measurable set $e\subset F$.

The two following facts are the main results\footnote{We note that
here we slightly changed the notation of the paper~\cite{Sue18} -- the new notation corresponds more fully
to the scalar approach based on the use of a~Riemann surface.}
of~\cite{Sue18}.

\begin{proposition}[\rm (see \cite{Sue18}, Theorem 1)]\label{prop1}
In the class $M_1^\circ(\FF)$, there exists a~unique measure
$\blambda=\lambda_{\FF}\in M_1^\circ(\FF)$ such that
\begin{equation}
J_\vv(\blambda)=\min_{\bmu\in M_1(\FF)}J_\vv(\bmu).
\label{18}
\end{equation}
The measure $\blambda$ is completely characterized by the following equilibrium conditions\footnote{The
equilibrium relations on the compact set~$\FF$ are stated somewhat
differently than in~\cite{Sue18}.
The equality $P^\blambda(\zz)+\vv(\zz)\equiv w_\FF$
everywhere on $S(\lambda_\FF)$ follows from the equality
$S(\lambda_\FF)=\FF$ (which is proved below) and the regularity of~$\FF$.}:
\begin{equation}
P^\blambda(\zz)+\vv(\zz)
\left\{\begin{matrix}
\,\leq w_{\FF}, & z\in S(\blambda),\\
\,\geq w_{\FF}, & z\in \FF\setminus S(\blambda).
\end{matrix}\right.
\label{19}
\end{equation}
Moreover, $P^\blambda(\zz)+\vv(\zz)\equiv w_\FF$ quasi-everywhere on $S(\blambda)$.
\end{proposition}

\begin{proposition}[\rm (see \cite{Sue18}, Theorem 2)]\label{prop2}
Let the functions $f_1$ and $f_2$ be given by~\eqref{1} and let $Q_{n,2}$
be the Hermite--Pad\'e polynomial defined by~\eqref{4}. Then
\begin{equation}
\frac1n\chi(Q_{n,2})\overset{*}\longrightarrow\lambda,\quad n\to\infty.
\label{20}
\end{equation}
\end{proposition}

The convergence in~\eqref{20} is understood in the sense of weak-$*$ convergence in the space of
measures, the measure $\lambda$ in~\eqref{20} is the measure
$\lambda=\pi_2(\lambda_{\FF})$ such that $\lambda(e)=\lambda_{\FF}(e^{(1)})$ for any
measurable set $e\subset F$.

Recall that the function $\pfi(\zz)=z+w$ is considered here as
a~function of a~point $\zz=(z,w)$ on the RS $\RS_2$. The
function~$\pfi$ is meromorphic on this RS (i.e., a~rational function
of~$z$ and~$w$, $\pfi\in\CC(z,w)$), has a~first-order pole at the point
$\zz=\infty^{(0)}$, and has a~first-order zero at the point
$\zz=\infty^{(1)}$. Therefore, its divisor $(\pfi)$ on $\RS_2$ is
$(\pfi)=-[\infty^{(0)}]+[\infty^{(1)}]$. The points $z=\pm1$ are
critical values of the canonical projection
$\pi_2\colon\RS_2\to\myh{\CC}$, $\pi_2(\zz)=z$. For the function~$\pfi$,
which is considered on the RS $\RS_2$, these points are regular. So,
the external field $\vv(\zz)=-\log|\pfi(\zz)|$ is harmonic on
$\RS_2\setminus\{\infty^{(0)},\infty^{(1)}\}$ and $P^\bmu(\zz)$ is
harmonic on~$\RS_2$ outside the set $F^{(0)}\cup
F^{(1)}\cup\infty^{(0)}\cup\infty^{(1)}$. Consequently, as distinct
from \eqref{5} and~\eqref{6}, the interval $E=[-1,1]$ is by no means
involved in both the definition of the function $P^\bmu(\zz)$ and in
the definition of the external field $\vv(\zz)$ (see a~slightly
modified statement of the equilibrium problem in~\cite{RaSu13}
and~\cite{MaRaSu16}, which nevertheless should be also considered as
a~vector problem, rather than a~scalar problem). So, by the above, the
equilibrium problem \eqref{19} can be naturally looked upon as a~scalar
equilibrium problem (but on a~two-sheeted RS and with harmonic external
field). We note the paper~\cite{Sue18b}, which deals with the case when
in the representation \eqref{1} for the function $f_2$ one considers,
as a~function~$h$, a~function holomorphic on~$E$ and having
in~$\myh{\CC}\setminus E$ a~finite number of branch points of arbitrary
character. In this setting, the branch points of the function~$h$ may
fail to be symmetric about the real line and the function~$h$ can now
assume complex values on $E$, and hence, as a~support of an appropriate
equilibrium measure, there naturally appears some $S$-compact set
$F^{*}$ (or, in a~different terminology, an $S$-curve; see~\cite{Rak12}
for more on this concept) instead of a~union of closed intervals of
the real line. The proof of the existence of such an $S$-compact set
$F^{*}$ is an involved problem (see first of all~\cite{Sta12}, and
also~\cite{BaStYa12} and~\cite{RaSu13}). Once the existence of
a~compact set $F^{*}$ is established, it will be used to define the
second compact set $E^{*}$ which also has the $S$-property and which is
a~natural replacement of the original interval $E=[-1,1]$. In some
cases (see~\cite{Sue18b}), the existence problem of a~compact
set~$F^{*}$ can be solved by an appeal to a~scalar problem in potential
theory of form \eqref{18}, which is posed on a~two-sheeted RS
(cf.~\cite{Apt08},~\cite{RaSu13}). This is a~certain advantage of the
scalar approach over the traditional vector approach.

The following fact will be required below. It is well known that
the solution $\vec{\lambda}=(\lambda_1,\lambda_2)$ of the vector
problem~\eqref{5}
exists and is unique (see \cite{Nik86},~\cite{NiSo88}, and
also \cite{ApBoYa17}, \cite{BaGeLo18}). Moreover, $S(\lambda_1)=E$ and
$S(\lambda_2)=F$; here
$\lambda_2=\beta_F(\lambda_1)$ is the balayage of the measure $\lambda_1$ from the domain
$\Omega=\myh{\CC}\setminus{F}$ onto $\partial\Omega=F$. Similarly,
$\lambda_1=\bigl(\beta_E(\lambda_2)+3\tau_E^{\hp}\bigr)/4$, where $\beta_E(\lambda_2)$
is the balayage of the measure~$\lambda_2$ from the domain $D=\myh{\CC}\setminus E$ onto $\partial
D=E$ and $\tau_E^{\hp}=dx/(\pi\sqrt{1-x^2})$ is the Chebyshev measure of the closed interval~$E$.

The purpose of the present paper is to show, first, that in~\eqref{19}
there is an equality on the whole of the compact set~$\FF$ (i.e., $S(\lambda_{\FF})=\FF$),
and second, to verify that the vector problem~\eqref{6} and the scalar problem~\eqref{19}
are equivalent. More precisely, the following result holds.

\begin{theorem}\label{th1}
Let $\lambda_{\FF}\in M_1^\circ(\FF)$ be a~(unique) unit measure
satisfying the equilibrium relations~\eqref{19} and let $\lambda_1\in M^\circ_1(E)$ be a~(unique)
unit measure satisfying the equilibrium relation~\eqref{6}. Then the following assertions hold:

1) $P^{\lambda_\FF}(\zz)+V(\zz)\equiv w_\FF$
 and $S(\lambda_{\FF})=\FF$ on the whole of the compact set~$\FF$;

2) the equilibrium problem~\eqref{6} and~\eqref{19} are equivalent, namely,
$\lambda=\beta_F(\lambda_1)$, where $\lambda=\pi_2(\lambda_\FF)$, and visa versa,
$\lambda_1=\bigl(\beta_E(\lambda)+3\tau_E^{\hp}\bigr)/4$,
where $\tau_E^{\hp}$ is the Chebyshev measure of the interval $E=[-1,1]$, $\beta_E(\lambda)$
is the balayage of the measure~$\lambda$ from the domain~$\Omega$ onto~$E$.
\end{theorem}

\begin{remark}\label{rem3}
Applying the operator $\dd^c$ to the function $P^\bmu(\zz)$, we get
$$
-\frac1{2\pi}\dd^c P^\bmu(\zz)=\bmu+\mu-\delta_{\infty^{(0)}}
-\delta_{\infty}.
$$
This shows that the function $P^\bmu$ is a~potential of the neutral charge
supported on the set $F^{(1)}\cup
F^{(0)}\cup\infty^{(0)}\cup\infty^{(1)}$ (see \cite{Chi18},~\cite{Chi19}).
\end{remark}

\begin{remark}\label{rem4}
It can be easily checked that the conclusions of Theorem~\ref{th1}
also hold in the more general setting when $F\subset\RR\setminus E$ is a~regular compact set.
\end{remark}

\section{Proof of Theorem~\ref{th1}}\label{s2}

\subsection{}\label{s2s1}
Let us show that $S(\lambda_\FF)=\FF$, and hence, in the equilibrium relations~\eqref{19}
the equality holds on the whole of the set~$\FF$.

Indeed, on the two-sheeted RS $\RS_2=\RS_2(w)$ with the above
Nuttall partition into (open) sheets $\RS_2^{(0)}$ and $\RS_2^{(1)}$ and
$\partial\RS_2^{(0)}=\partial\RS_2^{(1)}=:\bGamma$, where $\pi_2(\bGamma)=E$,
there is the involution operation~``$*$'' defined by $\zz=(z,w)\mapsto\zz^{*}:=(z,-w)$. This
operation swaps the sheets $\RS_2^{(0)}$ and $\RS_2^{(1)}$, but fixes the curve~$\bGamma$.

Given $\zz\in\RS_2^{(1)}$, we set
\begin{equation}
v(\zz):=P^{\lambda_{\FF}}(\zz)+\vv(\zz)-P^{\lambda_{\FF}}(\zz^*)-\vv(\zz^*).
\label{54}
\end{equation}
By the identity $\pfi(\zz)\pfi(\zz^*)\equiv1$ for $\zz\in\RS_2^{(1)}$, we have
\begin{align}
v(\zz)&=
\int_{\FF}\log\frac{\left|1-1/\bigl(\pfi(\zz)\pfi(t^{(1)})\bigr)\right|}
{|1-\pfi(\zz)/\pfi(t^{(1)})|}
\,d\lambda_{\FF}(t^{(1)})-2\log|\pfi(\zz)|\notag\\
&=\int_{\FF}\log\frac{|1-\pfi(\zz)\pfi(t^{(1)})|}
{|\pfi(\zz)-\pfi(t^{(1)})|}\,d\lambda_{\FF}(t^{(1)})-3\log|\pfi(\zz)|.
\label{55}
\end{align}
It follows that $v(\zz)$ is a~superharmonic function in the domain
$\RS_2^{(1)}\setminus\infty^{(1)}$, which extends continuously to~$\bGamma$,
because $\FF\cap\bGamma=\varnothing$ and since $v(\zz)\equiv0$ for $\zz\in\bGamma$.
We have $\pfi(\zz)\to0$
as $\zz\to\infty^{(1)}$, and hence from \eqref{55} we have $v(\zz)\to\infty$
as $\zz\to\infty^{(1)}$. Therefore, $v(\zz)>0$ for $\zz\in\RS_2^{(1)}$.
So, $P^{\lambda_{\FF}}(\zz)+\vv(\zz)>
P^{\lambda_{\FF}}(\zz^*)+\vv(\zz^*)$ for $\zz\in\FF$, because of
$\FF\cap\bGamma=\varnothing$. As a~result,
$P^{\lambda_{\FF}}(\zz^*)+\vv(\zz^*)<w_{\FF}$ for
$\zz\in S(\lambda_{\FF})=:S^{(1)}$. Hence
$P^{\lambda_{\FF}}(\zz)+\vv(\zz)<w_{\FF}$ for $\zz\in S^{(0)}$,
$S^{(0)}:=\{z^{(0)}\in\RS^{(0)}_2: z^{(1)}\in S^{(1)}\}$.
The function $u(\zz):=P^{\lambda_{\FF}}(\zz)+\vv(\zz)$ is harmonic in the domain $\mathfrak
D:=\RS_2\setminus(S^{(0)}\cup S^{(1)}\cup\infty^{(0)}\cup\infty^{(1)})$,
$u(\zz)\leq w_F$ for $\zz\in S^{(1)}$,
$u(\zz)<w_\FF$ for $\zz\in S^{(0)}$ and $u(\zz)=-3\log|z|+O(1)$ for
$\zz\to\infty^{(0)}$, $u(\zz)=O(1)$ for $\zz\to\infty^{(1)}$.
Applying the operator $\dd^c$ to the function~$u$, we see that
$$
-\frac1{2\pi}\dd^c u=\bmu+\mu-3\delta_{\infty^{(0)}}.
$$
Hence, the function $u$ is harmonic near the point $\zz=\infty^{(1)}$.
Therefore, $u(\zz)<w_F$ for $\zz\in\mathfrak D$. In particular, if
$\FF\setminus S^{(1)}\neq\varnothing$, then the inequality
$u(\zz)<w_\FF$ should be satisfied for $\zz\in\FF\setminus S^{(1)}$,
whereas by~\eqref{19} we have the reverse inequality $u(\zz)\geq
w_\FF$ for $\zz\in\FF\setminus S^{(1)}$. Hence $\FF\setminus
S^{(1)}=\varnothing$, $S(\lambda_{\FF})=\FF$ and
$P^{\lambda_{\FF}}(\zz)+\vv(\zz)\equiv w_{\FF}$ quasi-everywhere
on~$\FF$, $P^{\lambda_{\FF}}(\zz)+\vv(\zz)\leq w_{\FF}$ for
$\zz\in\FF$. Since $\FF$ is a~regular compact set, it follows that
$P^{\lambda_{\FF}}(\zz)+\vv(\zz)\equiv w_{\FF}$ everywhere on~$\FF$.

\subsection{}\label{s2s2}
Here we need some results from \cite{BuSu15}.
Even though, as one can easily check, the required relation is a~direct consequence of
\cite[formula (18)]{BuSu15} with $\theta=3$, for the sake of completeness
we provide its proof.

Let $\lambda_2=\beta_F(\lambda_1)$ be the balayage of the measure~$\lambda_1$ from the domain~$D$ onto~$F$. Since $F$ is a~regular compact set,
we have
\begin{equation}
U^{\lambda_2}(z)=U^{\lambda_1}(z)-G^{\lambda_1}_F(z)+\const,
\quad z\in\myh{\CC}
\label{41}
\end{equation}
(the value of the constant $\const$ is irrelevant here). Consider now the function
\begin{equation}
v(z):=3U^{\lambda_1}(z)+G^{\lambda_1}_F(z)+G^{\lambda_2}_E(z)+3g_E(z,\infty),
\quad z\in\myh{\CC},
\label{42}
\end{equation}
where, for an arbitrary measure $\nu$, $S(\nu)\subset\CC$,
$$
G^\nu_E(z):=\int g_E(z,\zeta)\,d\nu(\zeta),
\quad z\in D,
$$
is the Green potential of the measure $\nu$ and $g_E(z,\zeta)$ is the Green function for the domain~$D$
with singularity at the point $z=\zeta$.
By~\eqref{6} we have
\begin{equation}
v(z)=3U^{\lambda_1}(z)+G^{\lambda_1}_F(z)\equiv w_E=\const,
\quad z\in E.
\label{43}
\end{equation}
Since $\lambda_1$ is a unit measure,
$g_E(z,\infty)=\gamma_E-U^{\tau_E^{\hp}}(z)$, $\gamma_E=2\log{2}$ is the Robin constant
for $E$, $\tau_E^{\hp}$ is the Chebyshev measure for~$E$, it follows that the function~$v$
is harmonic in the domain $\myh{\CC}\setminus(E\cup F)$. Applying the
operator $\dd^c$ to both sides of~\eqref{42}, this establishes
\begin{equation}
-\frac1{2\pi}\dd^cv=3\lambda_1-3\delta_\infty+\lambda_1-\lambda_2+\lambda_2
-\beta_E(\lambda_2)-3\tau_E^{\hp}+3\delta_\infty,
\label{44}
\end{equation}
where $\beta_E(\lambda_2)$ is the balayage of the measure $\lambda_2$ from the domain
$\Omega$ onto~$E$. From \eqref{44} we see that
\begin{equation}
-\frac1{2\pi}\dd^cv=4\lambda_1-\beta_E(\lambda_2)-3\tau_E^{\hp},
\label{45}
\end{equation}
which shows that the function $v$ is harmonic already in the domain~$D$. So, by~\eqref{42} and~\eqref{45},
the function~$v$ is a~potential of the neutral charge with support on~$E $ and which is
identically constant on~$E$. The potential~$v$
is continuous on the support~$E$, and hence $v$~is continuous also in~$\myh{\CC}$. Thus, the function~$v$,
which is continuous on~$\myh{\CC}$, is harmonic on the domain~$D$ and is constant on
$\partial D=E$. Hence, $v$~is a~constant function.
Namely, $v(z)\equiv w_E$ for $z\in\myh{\CC}$. By~\eqref{41} we have
$U^{\lambda_1}(z)\equiv U^{\lambda_2}(z)+\const$ for $z\in F$, and hence
from \eqref{41} and~\eqref{42} we find that
\begin{equation}
v(z)=3U^{\lambda_2}(z)+G^{\lambda_2}_E(z)+3g_E(z,\infty)
\equiv\const,\quad z\in F
\label{46}
\end{equation}
(cf.~\cite[formula (18)]{BuSu15}).

Let us return back to the scalar equilibrium relation~\eqref{19} (recall
that we have already proved that $S(\lambda_\FF)=\FF$).

Since $\pfi(z^{(0)})\pfi(z^{(1)})\equiv1$, we have, for $\zz\in\RS^{(1)}_2$,
\begin{equation}
\frac{\left|1-1/\bigl(\pfi(\zz)\pfi(t^{(1)})\bigr)\right|}{|z-t|^2}
=\frac{|1-\pfi(z^{(0)})\pfi(t^{(0)})|}{|z-t|^2}.
\label{48}
\end{equation}
Therefore,
\begin{equation}
\log\frac{\left|1-1/\bigl(\pfi(\zz)\pfi(t^{(1)})\bigr)\right|}{|z-t|^2}
=\log
\frac{|1-\pfi(z^{(0)})\pfi(t^{(0)})|}{|z-t|^2}.
\label{53}
\end{equation}
Moreover, $\vv(\zz)=\log|\pfi(z^{(0)})|$ for $\zz\in\RS_2^{(1)}$.
So, by identifying an arbitrary measure $\bmu\in M_1(\FF)$ with the measure
$\pi_2(\bmu)\in M_1(F)$ and putting
$z^{(0)}=z$ and $t^{(0)}=t$, we have from \eqref{53} the following representation for
the function
$P^\bmu(\zz)$ with the external field $\vv(\zz)$ for $\zz\in\RS_2^{(1)}$ and $\bmu\in M_1(\FF)$ (cf.~\cite[formula (1.13)]{Sue18}):
\begin{equation}
P^\bmu(\zz)+\vv(\zz)
=\int_F\log\frac{|1-\pfi(z)\pfi(t)|}{|z-t|^2}\,d\mu(t)+\log|\pfi(z)|.
\label{49}
\end{equation}

Given an arbitrary measure $\mu\in M_1^\circ(F)$, consider the mixed
Green--logarithmic potential (cf.~\eqref{46}):
\begin{equation}
v_2(z;\mu):=3U^\mu(z)+G^{\mu}_E(z)+3g_E(z,\infty).
\label{50}
\end{equation}
Since
$E,F$ are compact subsets of $\RR$ and since $\pfi(z)$ is a~real-valued function for
$z\in\RR\setminus E$, we have, for the Green function $g_E(z,t)$ for $z,t\in\RR\setminus E$,
\begin{equation}
g_E(z,t)
=\log\frac{|1-\pfi(z)\myo{\pfi(t)}|}{|\pfi(z)-\pfi(t)|}
=\log\frac{|1-\pfi(z)\pfi(t)|}{|\pfi(z)-\pfi(t)|}.
\label{51}
\end{equation}

We now employ the following easily verified identity (see~\cite{Sue18}):
\begin{equation}
z-a\equiv -\frac{\bigl(\pfi(\zz)-\pfi(\maa)\bigr)
\bigl(1-\pfi(\zz)\pfi(\maa)\bigr)}
{2\pfi(\zz)\pfi(\maa)},
\quad z,a\in D.
\notag
\end{equation}
Using this relation and~\eqref{48}, we finally get from \eqref{51}
\begin{equation}
g_E(z,t)=\log\frac{|1-\pfi(z)\pfi(t)|^2}{2|z-t|\cdot|\pfi(z)\pfi(t)|}.
\label{52}
\end{equation}
Therefore,
\begin{align}
v_2(z;\mu)&=3\int\log\frac1{|z-t|}\,d\mu(t)
+\int\log g_E(z,t)\,d\mu(t)
+3\log|\pfi(z)|\notag\\
&=\int_F\log\frac{|1-\pfi(z)\pfi(t)|^2}{2|z-t|^4|\pfi(z)\pfi(t)|}\,d\mu(t)
+3\log|\pfi(z)|\notag\\
&=2\int_F\log\frac{|1-\pfi(z)\pfi(t)|}{|z-t|^2}\,d\mu(t)+2\log|\pfi(z)|
-\log2-\int_F\log|\pfi(t)|\,d\mu(t)\notag\\
&=2\(P^\mu(z)+\log|\pfi(z)|\)+\const.
\label{52}
\end{align}
Now from \eqref{42},~\eqref{43},~\eqref{49},~\eqref{50} and~\eqref{52} it follows that the
vector equilibrium problem~\eqref{6}
and the scalar equilibrium problem \eqref{19} are equivalent and that
$\lambda=\pi_2(\lambda_\FF)=\lambda_2$.

The function $v(z)$, as defined by~\eqref{42}, is identically constant, and hence,
applying the operator $\dd^c$ to both sides of~\eqref{42}, we get
\begin{equation}
0=-3\lambda_1+3\delta_\infty-\lambda_1+\lambda_2-\lambda_2
+\beta_E(\lambda_2)+\tau_E^{\hp}-3\delta_\infty
=\beta_E(\lambda_2)+3\tau_E^{\hp}-4\lambda_1.
\notag
\end{equation}
This implies the required representation
$\lambda_1=\dfrac14\beta_E(\lambda)+\dfrac34\tau_E^{\hp}$.

\end{document}